\documentclass[preprint,12pt]{elsarticle}
\usepackage{amsmath}
\usepackage{mathdots}
\usepackage{latexsym,amsfonts,amssymb}
\usepackage{lineno,hyperref}
\usepackage[top=3.2cm,bottom=2.2cm,left=2.5cm,right=2.5cm]{geometry}
\usepackage{natbib}
\usepackage{bbm}
\newtheorem{thm}{Theorem} 
\newtheorem{lem}[thm]{Lemma} 
\newdefinition{rmk}{Remark} 
\newproof{pf}{Proof}
\newproof{pot}{Proof of Theorem \ref{thm:LDP_PSDE}}
\newtheorem{dfn}[thm]{\quad \ \ Defination}
\newcommand{\kuo}[1]{\left( #1 \right)}
\newcommand{\zkuo}[1]{\left[ #1 \right]}
\newcommand{\dkuo}[1]{\left \{ #1 \right \}}
\newcommand{\e}{\mathbb{E}}
\newcommand{\p}{\mathbb{P}}
\newcommand{\G}{\mathcal{G}}
\newcommand{\abs}[1]{\left|#1\right|}
\newcommand{\norm}[1]{\|#1\|}
\makeatletter
\def\ps@pprintTitle{%
	\let\@oddhead\@empty
	\let\@evenhead\@empty
	\let\@oddfoot\@empty
	\let\@evenfoot\@oddfoot
}
\makeatother
\begin{document}

\begin{frontmatter}



\title{Large deviation for small noise path-dependent stochastic differential equations}

\author{Liu Xiangdong$^{\ast}$}
\cortext[mycorrespondingauthor]{Corresponding author}
\ead{tliuxd@jnu.edu.cn}
\author{Hong Shaopeng\corref{myauthor}}
\ead{hsp1999@stu2021.jnu.edu.cn}

\affiliation{organization={Department of Statistic and Data Science, School of Economics},
            addressline={Jinan University}, 
            city={Guangzhou},
            postcode={510632}, 
            state={Guangdong},
            country={PR China}}
\begin{abstract}
In this paper, we study the asymptotic behavior of randomly perturbed path-dependent stochastic differential equations with small parameter $\vartheta_{\varepsilon}$, when $\varepsilon \rightarrow 0$, $\vartheta_\varepsilon$ goes to $0$.
 When $\varepsilon \rightarrow 0$, we establish large deviation principle. The proof of the results relies on the weak convergence approach. As an application, we establish the large deviation for functionals of path-dependent SDEs in small time intervals.
\end{abstract}



\end{frontmatter}


\section{Introduction}
\label{sec:1}
This paper sheds new light on the asymptotic behaviour of the class of path-dependent stochastic  differential equations (PSDEs).
\begin{equation}\label{eq:1.1}
  X(t) = X_0 + \int_0^t  b(s,X_s)ds + \int_0^t \sigma(s,X_s)dW(s)\quad t \in [0,T]
\end{equation}

PSDEs have received increasing attentions by researchers which are much more involved than classical SDEs as the drift
and diffusion coefficients depending on path of solution. In a nutshell, this kind of equations plays an important role in characterising non-Markov  partial differential
equations  (PDEs for short). \citet{ekren2014viscosity} obtained the viscosity solutions of path-dependent semi-linear parabolic PDEs using backward PSDEs and Non-anticipative analysis \citep{dupireFunctionalItoCalculus2019,contFunctionalItoCalculus2013}, and subsequently extended the results to fully nonlinear forms of path-dependent PDEs \citep{Ekren2016ViscositySO}.

It is well known that the key point of large deviation principle (LDP for short) is to
show the probability property of rare events. Small noise LDP for SDEs has a long history. The pioneering work of \cite{freidlinRandomPerturbationsDynamical2012} considered rare events induced by Markov diffusions. Recently, an important contribution by \cite{budhiraja2019Analysis} was to use the weak convergence method to obtain a significant simplified approach. Their approach avoided proof exponential continuity and tightness estimates. 
Weakly convergent methods are widely used in proving large deviations of stochastic differential equations and stochastic partial differential equations, see \cite{jacquier2022large,Rckner2010LargeDF,xiong2018large} and references therein.

There have been some studies on large deviations of path-dependent SDEs. For example, \citet{gao2006LARGE} stuided such a problem via the sample path LDP method by Freidlin-Wentzell and show the LDP under (r,q)-capacity. And \citet{Ma2014LargeDF} based on PDEs method get the LDP of path-dependent SDEs. 
In this paper, we use a different line of argument, adapting the weak convergence approach of \citet{budhiraja2019Analysis} to the path-dependent case. 

Compared with the results mentioned above, the contribution of this paper is to study LDP when the coefficients of PSDEs are all depending on $\varepsilon$, i.e., the solutions of PSDEs possibly degenerate.
As an application, we establish the large deviation for functionals of PSDEs in small time intervals.

The paper is organized as follows. In Section \ref{sec:2}, we state the weak convergence method for the large deviation principle given in \citet{budhiraja2019Analysis}. We give the main theorem and prove it in Section \ref{sec:3}. Finally, in Section \ref{sec:4}, we show the large deviation principle for the functional of PSDEs in small time interval.

We end this section with some notations. We consider a fixed time horizon $T>0$, and denote $\mathbb{T}:=[0, T]$. Let $C([0, T] ; \mathbb{R}^d)$ be the Banach space of continuous functions $\psi:[0, T] \rightarrow \mathbb{R}^d$ equipped with the sup-norm $\|\psi\|:=\sup _{t \in[0, T]}|\psi(t)|$, $\mathcal{C}^1_0([0,T];\mathbb{R}^d)$ as the space of continuous functions on $[0,T]$ with initial value 0 and has first-order derivative, $\mathcal{C}^1_b([0,T];\mathbb{R}^d)$as the space of continuous functions on $[0,T]$ with initial value 0, has first-order derivative and has a bound. $L^2$ stands short for $L^2(\mathbb{T})$ and $\norm{\cdot}_2$ is the usual $L^2$ norm.

\section{Preliminaries}
\label{sec:2}
\subsection{Framework}
We consider small-noise convolution PSDEs
\begin{equation}\label{eq:2.1}
  X^\varepsilon(t) = X^\varepsilon_0 + \int_0^t b_\varepsilon(s,X^\varepsilon_s)ds + \vartheta_\varepsilon \int_0^t \sigma_\varepsilon(s,X^\varepsilon_s)dW(s)\quad t \in [0,T]
\end{equation}
taking values in $\mathbb{R}^d$ with $d\geq 1$, where $\varepsilon > 0$ and $\vartheta_\varepsilon >0 $ tends to zero as $\varepsilon $ goes to zero. For each $\varepsilon>0$, $X^\varepsilon_0 \in \mathbb{R}^d$, $b_\varepsilon : \mathbb{T}\times \mathcal{C}\kuo{\mathbb{T},\mathbb{R}^d}\rightarrow \mathbb{R}^d$, $\sigma_\varepsilon : \mathbb{T}\times \mathcal{C}\kuo{\mathbb{T},\mathbb{R}^d}\rightarrow \mathbb{R}^{d\times m}$are two product measurable maps that are non-anticipative  in the sense that they satisfy $b_\varepsilon(t,x) =  b_\varepsilon(t,x_t)$ and $\sigma_\varepsilon(t,x) = \sigma_\varepsilon(t,x_t)$ for all $t \in \mathbb{T}$  and each $x\in \mathcal{C}\kuo{\mathbb{T},\mathbb{R}^d}$, where $x_t$ denote the path $x$ stopped at time $t$. $W(s)$ is an m-dimensional Brownian motion on the filtered probability space $\kuo{\Omega,\mathcal{F},\dkuo{\mathcal{F}_t}_{t\in \mathbb{T}}, \mathbb{P}}$ satisfying the usual conditions.
We make following assumptions about the coefficients:
\begin{enumerate}
		\item [\textbf{A.1}] $X^\varepsilon_0$ converges to $x_0 \in \mathbb{R}^d$ as $\varepsilon$ tends to zero.
	\item [\textbf{A.2}] For all $\varepsilon > 0$ small enough, the coefficients $b_\varepsilon$ and $\sigma_\varepsilon$ are measurable maps on $\mathbb{T} \times \mathcal{C}\kuo{\mathbb{T}:\mathbb{R}^d}$ and converge pointwise to $b$ and $\sigma$ as $\varepsilon$ goes to zero. Moreover, $b(t,\cdot)$ and $\sigma(t,\cdot)$ are continuous on $\mathbb{R}^d$, uniformly in $t\in \mathbb{T}$.
	\item [\textbf{A.3}] For all $\varepsilon >0$ small enough, $b_\varepsilon$ and $\sigma_\varepsilon$ have linear growth uniformly in $\varepsilon$ and in $t \in \mathbb{T}$. For some $L>0$
	\begin{equation}
  |b_\varepsilon(t,\omega)| + |\sigma_\varepsilon(t,\omega) | \leq M \kuo{1+\sup_{s\leq t}|\omega(s)| + |t|}
\end{equation}
	\item [\textbf{A.4}] For all $\varepsilon > 0$ small enough, the coefficients $b_\varepsilon$ and $\sigma_\varepsilon$ are locally Lipschitz continuous. For any $R > 0$, there exists $L_R>0$ such that, for all 
\begin{equation}
  |b_\varepsilon(t,\omega) - b_\varepsilon(t,\omega^\prime) | + | \sigma_\varepsilon(t,\omega) - \sigma_\varepsilon(t,\omega^\prime) | \leq L_R (\sup_{s\leq t}|\omega(s) - \omega^\prime(s) |)
\end{equation}

\end{enumerate}
\subsection{Abstract sufficient conditions for large deviations}
\begin{dfn}[Large deviation \citep{demboLargeDeviationsTechniques2010}]{}
 A family $\left\{X^{\varepsilon}\right\}_{\varepsilon>0}$ of $\mathcal{E}$-valued random variable is said to satisfy the large deviation principle on $\mathcal{E}$, with the good rate function $I$ and with the speed function $\lambda(\varepsilon)$ which is a sequence of positive numbers tending to $+\infty$ as $\varepsilon \rightarrow 0$, if the following conditions hold:
 \begin{enumerate}
 	\item for each $M<\infty$, the level set $\{x \in \mathcal{E}: I(x) \leq M\}$ is a compact subset of $E$;
 	\item for each closed subset $F$ of $\mathcal{E}, \limsup _{\varepsilon \rightarrow 0} \frac{1}{\lambda(\varepsilon)} \log \mathbb{P}\left(X^{\varepsilon} \in F\right) \leq-\inf _{x \in F} I(x)$;
 	\item for each open subset $G$ of $\mathcal{E}, \liminf _{\varepsilon \rightarrow 0} \frac{1}{\lambda(\varepsilon)} \log \mathbb{P}\left(X^{\varepsilon} \in G\right) \geq-\inf _{x \in G} I(x)$.
 \end{enumerate}
 
\end{dfn}
We recall here several results from \citet{budhiraja2019Analysis} which gives an abstract framework of LDP.

Let $\mathcal{A}$ denote the class of real-valued $\left\{\mathcal{F}_t\right\}$-predictable processes $\nu$ belonging to $L^2$ a.s. For each $N$ the spaces of bounded deterministic and stochastic controls
$$
S_N:=\left\{\nu  \in L^2 ; \int_0^T|{\nu}(s)|^2 d s \leq N\right\} .
$$
$S_N$ is endowed with the weak topology induced from $L^2(\mathbb{T}\times \Omega)$. Define
$$
\mathcal{A}_N:=\left\{\nu \in \mathcal{A} ; \nu(s) \in S_N, \mathbb{P} \text {-a.s. }\right\} .
$$

\begin{thm}[\citet{budhiraja2019Analysis}]\label{thm:a_LDP}
For any $\varepsilon>0$, let $\mathcal{G}^{\varepsilon}$ be a measurable mapping from $C([0, T] ; \mathbb{R})$ into $E$. Suppose that $\left\{\mathcal{G}^{\varepsilon}\right\}_{\varepsilon>0}$ satisfies the following assumptions: there exists a measurable map $\mathcal{G}^0: C([0, T] ; \mathbb{R}) \longrightarrow \mathcal{E}$ such that
\begin{enumerate}[(a)]
	\item  for every $N<+\infty$ and any family $\left\{\nu^{\varepsilon} ; \varepsilon>0\right\} \subset \mathcal{A}_N$ satisfying that $\nu^{\varepsilon}$ converge in distribution as $S_N$-valued random elements to $\nu$ as $\varepsilon \rightarrow 0, \mathcal{G}^{\varepsilon}\left(W .+\frac{1}{\sqrt{\varepsilon}} \int_0^{\cdot} \nu^{\varepsilon}(s) d s\right)$ converges in distribution to $\mathcal{G}^0\left(\int_0^{\cdot} \nu(s) d s\right)$ as $\varepsilon \rightarrow 0$;
	\item for every $N<+\infty$, the set $\left\{\mathcal{G}^0\left(\int_0^\cdot \nu(s) d s\right) ; h \in S_N\right\}$ is a compact subset of $E$. \end{enumerate}
Then the family $\left\{\mathcal{G}^{\varepsilon}(W(\cdot))\right\}_{\varepsilon>0}$ satisfies a large deviation principle with the good rate function I given by
$$
I(g):=\inf _{\left\{\nu \in \mathcal{H} ; g=\mathcal{G}^0\left(\int_0^\cdot \nu (s) d s\right)\right\}}\left\{\frac{1}{2} \int_0^T|\nu(s)|^2 d s\right\} \quad \text { for } g \in \mathcal{E},
$$
with the convention $\inf \emptyset=\infty$.
\end{thm}
\section{Main Result and Proof}\label{sec:3}

If \textbf{A.5} hold, define the functional $\mathcal{G}$ as the Borel-measurable map associating the multidimensional Brownian motion $W$ to the solution of the path dependent stochastic differential systems \eqref{eq:2.1}, that is: $\G^\varepsilon \kuo{W} = X^\varepsilon(t) $. For any control $\nu \in \mathcal{A}_N$, $N>0$ and any $\varepsilon >0$, the process $\widetilde{W} = W + \vartheta_\varepsilon ^{-1}\int_0^{\cdot}\nu(s)ds$ is a $\widetilde{\p}-$Brownian motion by Girsanov's theorem, where
\begin{equation}
  \frac{d\widetilde{\p}}{d\p} :=\exp\dkuo{-\frac{1}{\vartheta_\varepsilon}\sum_{i=1}^m\int_0^T \nu^{(i)}(s)dW^{(i)}(s) - \frac{1}{2\vartheta_\varepsilon^2}\int_0^T|v(s)|^2ds} .
\end{equation}

Hence the shifted version $X^{\varepsilon, v}:=\mathcal{G}^{\varepsilon}(\tilde{W})$ appearing in Theorem \ref{thm:a_LDP} (1) is the strong unique solution of \eqref{eq:2.1} under $\widetilde{\mathbb{P}}$, with $X^{\varepsilon}$ and $W$ replaced by $X^{\varepsilon, v}$ and $\widetilde{W}$. Because $\mathbb{P}$ and $\widetilde{\mathbb{P}}$ are equivalent, $X^{\varepsilon, v}$ is also the unique strong solution, under $\mathbb{P}$, of the controlled equation
\begin{equation}\label{eq:3.1}
  X^{\varepsilon, v}(t)=X_0^{\varepsilon}+\int_0^t \left[b_{\varepsilon}\left(s, X_s^{\varepsilon, v}\right)+\sigma_{\varepsilon}\left(s, X_s^{\varepsilon, v}\right) v(s)\right] \mathrm{d} s+\vartheta_\varepsilon\int_0^t \sigma_{\varepsilon}\left(s, X_s^{\varepsilon, v}\right) \mathrm{d} W(s)
\end{equation}
Taking $\varepsilon \rightarrow 0$, the system \eqref{eq:3.1} reduces to the deterministic path dependent ODE 
\begin{equation}\label{eq:3.2}
  \phi(t) = x_0 + \int_0^t\zkuo{b(s,\phi_s)+\sigma(s,\phi_s)\nu(s)}ds.
\end{equation}

\begin{thm}{}\label{thm:LDP_PSDE}
Under \textbf{A.1}-\textbf{A.4}, the family $\dkuo{X^\varepsilon}_{\varepsilon>0}$, unique solution of \eqref{eq:2.1}, satisfies a Large Deviations Principle with rare function $I$ and speed $\vartheta_\varepsilon^{-2}$, where $\mathcal{G}^0$ is the solution of \eqref{eq:3.2} 
\end{thm}
\begin{rmk}
Theorem \ref{thm:LDP_PSDE} generalizes the results in \citet{chiarini2014Large}. When the coefficients $b_\varepsilon$ and $\sigma_\varepsilon $ do not depend on the path of the process $X_\varepsilon$, Theorem  \ref{thm:LDP_PSDE} and 
 \citet[Theorem 3]{chiarini2014Large} are equivalent.

\end{rmk}

We first show the unique solution of \eqref{eq:3.2} and a uniform estimation.
\begin{lem}
	Under \textbf{A.1}-\textbf{A.6}, given any $\nu \in L^2$, there is a unique solution $\phi \in \mathcal{C}\kuo{[0,T];\mathbb{R}^n}$ of \eqref{eq:3.2}. Moreover, for $\phi$, we have the growth estimate 
	\begin{equation}\label{eq:phi_bound}
  \sup_{0\leq s\leq t}\abs{\phi(s)}\leq (3|x_0|^2 + 9M^2t (t + \norm{\nu}^2)+ 3M^2t^3(t + \norm{\nu}^2))e^{9M^2(t + \norm{\nu}^2)}
\end{equation}
\end{lem}
\begin{pf}
	Let $\phi, \varphi \in \mathcal{C}\kuo{[0,T];\mathbb{R}^d}$ be solutions of \eqref{eq:3.2}. We have 
	\begin{equation}
  \begin{aligned}
  	\abs{\phi(t) - \varphi(t)}\leq \int_0^t \abs{b(s,\phi_s) - b(s,\varphi_s) }ds + \int_0^t \abs{\sigma(s,\phi_s) - \sigma(s,\varphi_s)}\abs{\nu}ds
  \end{aligned}
\end{equation}
 By assumption \textbf{A.4}, we have for large enough $R>0$
 $$
 \abs{\phi(t)-\varphi(t)}^2\leq 2L^2_R\kuo{T+\norm{\nu}^2}\int_0^t\sup_{0\leq u \leq s}\abs{\phi(u)-\phi(s)}^2ds
 $$
 Gronwall's inequality now entails that $\norm{\phi(t)-\varphi(t)} = 0$, which yields uniqueness.
 
 By using assumption \textbf{A.3}, we can get that 
 \begin{equation}
  \begin{aligned}
  	\abs{\phi(t)}^2&\leq 3|x_0| + 3t\int_0^t\abs{b(s,\phi_s)^2}ds+3\kuo{\int_0^t\abs{\sigma(s,\phi_s)}|\nu|ds}^2\\
  	& \leq 3|x_0|^2 + 9 M^2\kuo{t+\norm{\nu}^2}\int_0^t\kuo{1+\sup_{0\leq u\leq s}|\phi(u)|^2+|s|^2}ds\\
  	& \leq 3|x_0|^2 + 9M^2t (t + \norm{\nu}^2)+ 3M^2t^3(t + \norm{\nu}^2) + 9M^2(t + \norm{\nu}^2)\int_0^t\sup_{0\leq u\leq s}|\phi(u)|^2du 
  	  	\end{aligned}
\end{equation}
By Gronwalls' inequality, we can deduce
$$
\sup_{0\leq s\leq t}\abs{\phi(s)}\leq (3|x_0|^2 + 9M^2t (t + \norm{\nu}^2)+ 3M^2t^3(t + \norm{\nu}^2))e^{9M^2(t + \norm{\nu}^2)}
$$

\end{pf}

We need some technical preliminary results.
\begin{lem}\label{lem:1}
	Under \textbf{A.1}-\textbf{A.6}, for all $p\geq 2$, $N >0$, $\nu \in \mathcal{A}_N$ and $\varepsilon >0$ small enough, there exists a constant $c>0$ independent of $\varepsilon$, $\nu$, $t$ such that 
	\begin{equation}
  \e\zkuo{\sup_{t \in \mathbb{T}}|X^{\varepsilon, \nu}(t)|^p}\leq c
\end{equation}
\end{lem}
\begin{pf}
	Let us fix $p\geq 2$, $N >0$, $\nu \in \mathcal{A}_N$ and $t \in \mathbb{T}$. Let $\tau_n$ be the stopping time defined by 
	$$
	\tau_n = \inf\dkuo{t\geq 0: |X^{\varepsilon,\nu}(t)|\geq n}\wedge T
	$$
	We write $b^n_s := b_\varepsilon(s, X^{\varepsilon,\nu}_s \mathbbm{1}_{\dkuo{s\leq \tau_n}})$  and $\sigma^n_s := \sigma_{\varepsilon}\kuo{s,X^{\varepsilon,\nu}_s \mathbbm{1}_{\dkuo{s\leq \tau_n}}}$. 
	
	We fix $n \in \mathbb{N}$ and observe that, almost surely: 
	\begin{equation}
  \begin{aligned}
  	\e\zkuo{\norm{X^{\varepsilon,\nu}(t)\mathbbm{1}_{t \leq \tau_n}}^p}&\leq 4^{p-1}\abs{X^\varepsilon_0}^p + 4^{p-1}\e\dkuo{\zkuo{\int_0^tb_s^nds}^p}+4^{p-1}\e\dkuo{\zkuo{\int_0^t\sigma_s^n\nu(s)ds}^p}\\
  	& + 4^{p-1}\vartheta_\varepsilon^p\e\dkuo{\zkuo{\int_0^t\sigma^n_sdW(s)}^p}\\
  	&=: 4^{p-1}\kuo{\abs{X_0^\varepsilon}^p + I_1+I_2+I_3}
  \end{aligned}
\end{equation}

For $\varepsilon$ small enough we can bound $\abs{X^\varepsilon_0}$ by $2\abs{X_0}$ and $\vartheta_\varepsilon$ by $1$. Using Holder's and Jensen's inequalities, we obtain the following estimates almost surely:
\begin{equation}
\begin{aligned}
	I_1&\leq \e\dkuo{\zkuo{\int_0^t(b_n^s)^pds}}
\end{aligned}
\end{equation}
and
\begin{equation}
  \begin{aligned}
  	I_2 & \leq N^{\frac{p}{2}}\e\dkuo{\zkuo{\int_0^t(\sigma_s^n)^2ds}^{\frac{p}{2}}}\leq N^{\frac{p}{2}}\e\dkuo{\int_0^t \kuo{\sigma_s^n}^pds}
  \end{aligned}
\end{equation}
By Burkholder-Davis-Gundy (B-D-G) inequality, there exists $C_p >0$ such that 
\begin{equation}
  \e\zkuo{I_3}\leq C_p \e\dkuo{\int_0^t(\sigma_s^n)^Pds}
\end{equation}
From the linear growth condition on $b_\varepsilon$ and $\sigma_\varepsilon$ we deduce that there exists $C_1>0$ independent of $\varepsilon$, $\nu$, $n$ and $t$ such that for all $n\in \mathbb{N}$
\begin{equation}
  \e\zkuo{\norm{X^{\varepsilon, \nu}(t)\mathbbm{1}_{t \leq \tau_n}}^p}\leq C_1 + C_1 \int_0^t \e\zkuo{\norm{X^{\varepsilon, \nu}(s)\mathbbm{1}_{s \leq \tau_n}}^p}ds
\end{equation}

Taking $n$ goes to infinity and using Gronwall's lemma, we prove this bound.
\end{pf}
\begin{lem}\label{lem:2}
$\dkuo{X^{\varepsilon,\nu^\varepsilon}}_{\varepsilon>0}$ is tightness	
\end{lem}
\begin{pf}
	In view of the Kolmogorov tightness criterion, it suffices to show that there exist strictly positive constants $\alpha$, $\beta$ and $\gamma$ such that for all $t$, $s\in [0,T]$,
	$$
	\sup_{\nu\in \mathcal{S}_N}\e\zkuo{\abs{X^{\varepsilon,\nu^\varepsilon}(t) - X^{\varepsilon,\nu^\varepsilon}(s)}^\alpha}\leq \beta\abs{t-s}^\gamma
	$$
	Without loss of generality, let $s<t$. We will write $b^n_s := b_\varepsilon(s, X^{\varepsilon,\nu^\varepsilon}_s \mathbbm{1}_{\dkuo{s\leq \tau_n}})$  and $\sigma^n_s := \sigma_{\varepsilon}\kuo{s,X^{\varepsilon,\nu^\varepsilon}_s \mathbbm{1}_{\dkuo{s\leq \tau_n}}}$.
	\begin{equation}
  \begin{aligned}
  	\e\zkuo{\abs{X^{\varepsilon,\nu^\varepsilon}(t) - X^{\varepsilon,\nu^\varepsilon}(s)}^\alpha}&\leq 3^{p-1}(t-s)^{p-1}\e\dkuo{\int_s^t\abs{b^n_u}^pdu}+3^{p-1}\e\zkuo{\kuo{\int_s^t\abs{\sigma_u^n}\abs{\nu(u)}du}^p}\\
  	&+ 3^{p-1}\vartheta_\varepsilon^P\e\zkuo{\abs{\int_s^t\sigma_u^ndW(u)}^p}\\
  	&\leq 3^{p-1}(t-s)^{p-1}\e\dkuo{\int_s^t\abs{b^n_u}^pdu} +3^{p-1}N^{\frac{p}{2}}\e\zkuo{\int_s^t\abs{\sigma^n_u}^pdu}\\
  	&+3^{p-1}C_p\e\zkuo{\int_s^t|\sigma^n_u|^{\frac{p}{2}}du}\\
  	&\leq 3^{p-1}(t-s)^{p-1}\e\dkuo{\int_s^t\abs{b^n_u}^pdu} + 3^{p-1}N^{\frac{p}{2}}\e\zkuo{\int_0^t |\sigma^n_u|^{\frac{p}{2}}du}\\
  	&+3^{p-1}C_p\e\zkuo{\int_s^t|\sigma^n_u|^{\frac{p}{2}}du}
  \end{aligned}
\end{equation}
From the linear growth condition on $b_\varepsilon$ and $\sigma_\varepsilon$ we deduce that there exists sufficiently large $\beta$ and let $\alpha = p$, $\gamma = p-1$. Then the hypotheses of Kolmogorov's criterion are therefore satisfied.
\end{pf}
\begin{lem}\label{lem:3}
	For any positive $N < \infty$, the set 
$$
K_N : = \dkuo{\mathcal{G}^0\kuo{\int_0^\cdot \nu(s) ds , \nu \in \mathcal{S}_N}}
$$
is a compact set in $\mathcal{C}\kuo{[0,T];\mathbb{R}^n}$
\end{lem}
\begin{pf}
	We first prove $\mathcal{G}^0$ is a continuity map form $\mathcal{S}_N$ to $\mathcal{C}\kuo{[0,T];\mathbb{R}}$, then for any positive $N< \infty$, $\mathcal{S}_N$ is compact set in weak topology. Since $\mathcal{G}^0$ is continuity map, we can show $K_N$ is a compact set in $\mathcal{C}\kuo{[0,T];\mathbb{R}}$.
	
	Taking $\dkuo{\nu^n(s)}\in \mathcal{S}_N$, $\nu^n \rightarrow \nu$ weakly, let $\varphi^n = \mathcal{G}^0(\nu^n)$, $\varphi = \mathcal{G}^0(\nu)$. Then, for $t \in [0,T]$, 
	\begin{equation}
  \begin{aligned}
  	\varphi^n(t) - \varphi(t) &= \int_0^t \kuo{b(s,\varphi^n) - b(s,\varphi)}ds + \int_0^t(\sigma(s,\varphi^n) - \sigma(s,\varphi))\nu^n(s)ds \\
  	&+ \int_0^t\sigma(s,\varphi_s)(\nu^n(s)-\nu(s))ds
  \end{aligned}
\end{equation}
Since $\norm{\nu^n}\leq N$, it follows from \eqref{eq:phi_bound} that $R := \sup_{n\in \mathbb{N}} \norm{\varphi}\vee \norm{\varphi^n}$ is finite. Therefore, using assumption \textbf{A.4}, 
\begin{equation}
  \begin{aligned}
  	\sup_{0\leq s\leq t}\abs{\varphi^n(s) - \varphi(s)}&\leq L_R\int_0^t\sup_{0\leq u \leq s}\abs{\varphi^n(u) - \varphi(u)}ds +L_R\int_0^t\sup_{0\leq u \leq s}\abs{\varphi^n(u) - \varphi(u)}\nu^n(s) ds\\
  	&+ \sup_{0\leq u\leq T}\abs{\int_0^u \sigma(s,\varphi_s)\kuo{\nu^n(s)-\nu(s)}ds}
  \end{aligned}
\end{equation}
Let $\Delta^n_\sigma = \sup_{0\leq u\leq T}\abs{\int_0^u \sigma(s,\varphi_s)\kuo{\nu^n(s)-\nu(s)}ds}$. By Holder's inequality and since $\norm{\nu^n}^2 \leq N$ for all $n \in \mathbb{N}$, it follows that 
$$
\sup_{0\leq s\leq t}\abs{\varphi^n(s) - \varphi(s)} \leq 3L_R^2(t +N) \int_0^t \sup_{0\leq u\leq s}\abs{\varphi^n(u) - \varphi(u)}^2ds + 3 (\Delta^n_\sigma )^2
$$

By Gronwall's lemma, we can deduce that 
$$
\mathcal{G}^0(\nu^n)-\mathcal{G}^0(\nu) = \sup_{0\leq t\leq T}|\varphi^n(t) - \varphi(t)|^2 \leq 3\kuo{\Delta^n_\sigma}^2e^{3L^2T(T+N)}
$$
In order to establish continuity of $\mathcal{G}^0$ on $\mathcal{S}_N$, it remains to check that $\Delta^n_\sigma$ goes to $0$ as $n\rightarrow \infty$. By \textbf{A.3}, it follows that $\sigma(\cdot,\varphi)\nu^n$ converges weakly to $\sigma\kuo{\cdot,\varphi}\nu$ in $L^2$. Moreover, the family $\dkuo{\sigma(\cdot,\varphi)\nu^n}_{n\in \mathbb{N}}$ is bounded in $L^2$ with respect to the $L^2-$norm. Hence, 
$$
\int_0^t \sigma(s,\varphi_s)\nu^n(s)ds \rightarrow \int_0^t \sigma\kuo{s,\varphi}\nu(s)ds \quad\text{as }n\rightarrow \infty
$$ 
which implies that $\Delta^n_\sigma \rightarrow 0$ as $n\rightarrow \infty$.
\end{pf}
\begin{lem}\label{lem:4}
	Under \textbf{A.1}-\textbf{A.6}, for every $N<+\infty$ and any family $\dkuo{\nu^\varepsilon}_{\varepsilon>0} \in \mathcal{A}_N$ satisfying that $\nu^\varepsilon $ converge in distribution as $\mathcal{S}_N-$valued random elements to $\nu$ as $\varepsilon \rightarrow 0$, $\mathcal{G}^\varepsilon\kuo{W_{\cdot}+ \frac{1}{\vartheta_\varepsilon}\int_0^\cdot \nu^\varepsilon(s)ds}$ converges in distribution to $\mathcal{G}^0\kuo{\int_0^\cdot \nu(s) ds}$ as $\varepsilon \rightarrow 0$.
\end{lem}
\begin{pf}
	By Skorohod representation theorem we can work with almost sure convergence for the purpose of identifying the limit. We follow the technique in \citet{chiarini2014Large}.

	For $t \in [0,T]$, define $\Phi_t:\mathcal{S}_N\times \mathcal{C}\kuo{[0,T]; R^n}$ as 
	$$
	\Phi_t(\omega,f) := \abs{\omega(t) - x_0 - \int_0^t b(s,\omega_s)ds - \int_0^t \sigma(s,\omega_s)f(s)ds}\wedge 1
	$$
	$\Phi_t$ is bounded and we show that it is also continuous. Let $\omega^n\rightarrow \omega$ in $\mathcal{C}\kuo{[0,T]; R^n}$ and $f^n \rightarrow f$ in $\mathcal{S}_N$ with respect to the weak topology. \textbf{A.2} implies the existence of continuous moduli of continuity $\rho_b$ and $\rho_\sigma$ for both coefficients such that $\abs{b(t,\varphi_t) - b(t,\phi_t)}\leq \rho_b\kuo{\norm{\varphi - \phi}}$ and $\abs{\sigma(t,\varphi_t) - \sigma(t,\phi_t)}\leq \rho_\sigma\kuo{\norm{\varphi - \phi}}$. Using Holder's inequality we find that 
	\begin{equation}
  \begin{aligned}
  	\abs{\Phi_t(\omega^n,f^n) - \Phi_t(\omega^n,f)}&\leq \abs{\omega^n(t) - \omega(t)} + \int_0^t\abs{b(s,\omega^n_s) - b(s,\omega_s)}ds \\
  	&+ \int_0^t \abs{\sigma(s,\omega^n_s) - \sigma(s,\omega_s)}|f^n(s)|ds+\abs{\int_0^t\sigma(s,\omega_s)(f^n(s)-f(s))ds}\\
  	&\leq \norm{\omega^n-\omega}+ T \rho_b\kuo{\norm{\omega^n-\omega}}+ \sqrt{NT}\rho_\sigma\kuo{\norm{\omega^n-\omega}}\\
  	&+\norm{\sigma(\cdot,\omega)}\abs{\int_0^t\kuo{f^n(s) -f(s)}ds}
  \end{aligned}
\end{equation}
Since $f^n$ tends to $f$ weakly in $L^2$ then the last integral converges to zero as $n$ goes to infinity. Moreover $\lim\limits_{n\uparrow \infty}\norm{\omega^n-\omega} = 0$, which proves that $\Phi_t$ is continuous, and therefore
$$
\lim\limits_{n\uparrow \infty} \e\zkuo{\Phi_t(X^n,\nu^n)} = \e\zkuo{\Phi_t(X,\nu)}
$$

Define $b^R_\varepsilon:\zkuo{0,T}\times \mathcal{C}([0,T]; \mathbb{R}^d)$ and $\sigma^R_\varepsilon:\zkuo{0,T}\times \mathcal{C}([0,T]; \mathbb{R}^d)$ by 
$$
b^R_\varepsilon (s,\omega_s) = \left\{\begin{matrix}b_\varepsilon (s,\omega_s)\quad \text{if} \|\omega \|\leq R
 \\b_\varepsilon (s,\frac{R}{\|\omega \|}\omega _s )\quad  \text{otherwise}
\end{matrix}\right.
\quad 
\sigma^R_\varepsilon (s,\omega_s) = \left\{\begin{matrix}\sigma_\varepsilon (s,\omega_s)\quad \text{if} \|\omega \|\leq R
 \\\sigma_\varepsilon (s,\frac{R}{\|\omega \|}\omega _s )\quad  \text{otherwise}
\end{matrix}\right.
$$
It is clear that the function $b^R_\varepsilon$ and $\sigma^R_\varepsilon$ are globally Lipschitz and bounded. By assumption \textbf{A.2}, $b^R_\varepsilon \rightarrow b^R$ and $\sigma^R_\varepsilon \rightarrow \sigma^R$ uniformly on $[0,T]\times \mathcal{C}\kuo{[0,T];\mathbb{R}^n}$. In analogy with $\Phi_t$, set 
	$$
	\Phi^R_t(\omega,f) := \abs{\omega(t) - x_0 - \int_0^t b^R(s,\omega_s)ds - \int_0^t \sigma^R(s,\omega_s)f(s)ds}\wedge 1
	$$
	
	Consider the family $\dkuo{X^{R,\varepsilon,\nu}}$ of solutions to the PSDE
	$$
	 X^{R,\varepsilon, v}(t)=X_0^{\varepsilon}+\int_0^t \left[b^R_{\varepsilon}\left(s, X_s^{R,\varepsilon, \nu}\right)+\sigma^R_{R,\varepsilon}\left(s, X_s^{R,\varepsilon, \nu}\right) v(s)\right] \mathrm{d} s+\vartheta_\varepsilon\int_0^t \sigma^R_{R,\varepsilon}\left(s, X_s^{R,\varepsilon, \nu}\right) \mathrm{d} W(s)
	$$
	We will show 
	$$
	\lim\limits_{\varepsilon \rightarrow 0 }\e\zkuo{\Phi_t^R\kuo{X^{R,\varepsilon,\nu^\varepsilon},\nu}} = 0
	$$
	\begin{equation}
  \begin{aligned}
  	\e\zkuo{\Phi_t^R\kuo{X^{R,\varepsilon,\nu^\varepsilon},\nu}}&\leq\abs{X^\varepsilon_0 - X_0}+\e\zkuo{\int_0^t\abs{b^R_\varepsilon(s,X^{R,\varepsilon,\nu^\varepsilon})-b^R(s,X^{R,\varepsilon,\nu^\varepsilon})}ds}\\
  	&+\e\zkuo{\int_0^t\abs{\sigma^R_\varepsilon(s,X^{R,\varepsilon,\nu^\varepsilon})-\sigma^R(s,X^{R,\varepsilon,\nu^\varepsilon})}|\nu(s)|ds}\\
   &+\vartheta_\varepsilon \e\zkuo{\abs{\int_0^t\sigma^R_\varepsilon\kuo{s,X^{R,\varepsilon,\nu^\varepsilon}}dW(s)}}\\
  	&\leq \abs{X^\varepsilon_0 - X_0}+ t \norm{b^R_\varepsilon - b^R}+\norm{\sigma^R_\varepsilon - \sigma^R}\e\zkuo{\int_0^T|\nu(s)|ds}\\
   &+\vartheta_\varepsilon\sqrt{\int_0^t\e\zkuo{\sigma_\varepsilon^R(s,X^{R,\varepsilon,\nu^\varepsilon})^2}ds}
  \end{aligned}
\end{equation}

The last term in the above display tends to 0 since 
\begin{equation}
  \begin{aligned}
  	\sup_{\nu \in \mathcal{S}_N}\int_0^t\e\zkuo{\sigma_\varepsilon^R(s,X^{R,\varepsilon,\nu^\varepsilon})^2}ds&\leq 2 \sup_{\nu \in \mathcal{S}_N} \int_0^T\e\zkuo{\abs{\sigma^R(s,X^{R,\varepsilon,\nu^\varepsilon})}^2}ds\\
  	&+ 2  \sup_{\nu \in \mathcal{S}_N} \int_0^T\e\zkuo{\abs{\sigma^R_\varepsilon(s,X^{R,\varepsilon,\nu^\varepsilon})-\sigma^R(s,X^{R,\varepsilon,\nu^\varepsilon})}^2}ds\\
  	&\leq 2T \sup_{\nu \in \mathcal{S}_N}\norm{\sigma_{\varepsilon}^R - \sigma^R}^2 + 2 \sup_{\nu \in \mathcal{S}_N}\int_0^T\e\zkuo{\abs{\sigma^R(S,X^{R,\varepsilon,\nu^\varepsilon})}^2}ds\\
   &<\infty
  \end{aligned}
\end{equation}
Then, we have 
	$$
	\lim\limits_{n\rightarrow \infty}\e\zkuo{\Phi_t^R\kuo{X^{R,\varepsilon,\nu^\varepsilon},\nu}} = 0
	$$
	For $R>0$ and $\nu \in \mathcal{S}_N$, let $\tau^R_n$ is a stopping time defined by 
	$$
	\tau^R_n = \inf\dkuo{t\geq 0: X^{\varepsilon,\nu}(t) \geq R}
	$$
	We have 
	$$
	\p\kuo{X^{R,\varepsilon,\nu^\varepsilon}(t) = X^{\varepsilon,\nu^\varepsilon}(t)\mathbbm{1}_{t\leq \tau_n}} = 1
	$$
	It follows that 
	\begin{equation}\label{eq:Phi}
  \begin{aligned}
  	\e\zkuo{\Phi_t(X^{\varepsilon ,\nu},\nu)} &= \e\zkuo{\mathbbm{1}_{t<\tau_n}\Phi_t(X^{\varepsilon,\nu},\nu)}+\e\zkuo{\mathbbm{1}_{t\geq \tau_n}\Phi_t(X^{\varepsilon,\nu},\nu)}\\
  	&\leq \e\zkuo{\Phi^R_t(X^{R,\varepsilon,\nu^\varepsilon},\nu)}+\p\kuo{t \geq \tau_R^n}
  \end{aligned}
\end{equation}
For all $\nu \in \mathcal{S}_N$, by Markov's inequality we have 
$$
\p\kuo{t\geq \tau^R_n} = \p\kuo{\sup_{0\leq s\leq t}|X^{R,\varepsilon,\nu^\varepsilon}(s)|\geq R}\leq\frac{c}{R^2}
$$
Taking upper limits on both side of \eqref{eq:Phi}, we obtain
$$
\limsup_{\varepsilon \rightarrow 0 }\e\zkuo{\Phi_t(X^{\varepsilon,\nu},\nu)}\leq\limsup_{n\rightarrow \infty}\p\kuo{t\geq \tau^R_n}\leq \frac{c}{R}
$$
Since $R >0$ has been chosen arbitrarily, it follows that 
$$
\lim\limits_{\varepsilon \rightarrow 0 }\e\zkuo{\Phi_t(X^{\varepsilon,\nu^\varepsilon},\nu^\varepsilon )} = 0
$$
\end{pf}
\begin{pot}
According to Theorem \ref{thm:a_LDP}, combined with Lemma \ref{lem:1}, \ref{lem:2}, \ref{lem:3} and \ref{lem:4} , it can be seen that Theorem \ref{thm:LDP_PSDE} holds.	
\end{pot}

\section{Application: Small time large deviation principle for path-dependent stochastic differential equation}\label{sec:4}
In this section, we study the LDP for functional of PSDEs in small time intervals: $\dkuo{X(t), t\in \mathbb{T}}$ as $t \rightarrow 0$, where 
$$
X(t) = x_0 + \int_0^t b(s,X_s)ds + \int_0^t \sigma(s,X_s)dW(s)
$$
We rescale the small time problem to a small perturbation problem.

\begin{equation}\label{eq:4.2}
	\begin{aligned}
		X({\varepsilon t}) &= x_0 + \int_0^{\varepsilon t} b(s,X_s)ds + \int_0^{\varepsilon t}\sigma(s,X_s)dW(s)\\
		& = x_0 + \varepsilon \int_0^t b(\varepsilon t, X_{\varepsilon t})ds + \sqrt{\varepsilon}\int_0^{\varepsilon t}\sigma(\varepsilon t, X_{\varepsilon t})d\widehat{W}(s)\\
	\end{aligned}
\end{equation}

Where, $\widehat{W}(s) = \frac{1}{\sqrt{\varepsilon}}W(\varepsilon s )$. Let $U(t) = X(\varepsilon t)$, by \eqref{eq:4.2}, we have 
$$
U(t) = x_0 = \int_0^t b_\varepsilon(s, U_s)ds + \sqrt{\varepsilon}\int_0^t \sigma(s,U_s)d\widehat{W}(s)
$$

Now, we can use Theorem \ref{thm:LDP_PSDE} to obtain the LDP for small time PSDEs.
\begin{thm}{}{}
The process $X(\varepsilon t )$ satisfies LDP as $\varepsilon \rightarrow 0$ with rate function $J^1$ and speed $\varepsilon $.
$$
J(g) = \inf_{\nu \in L^2; g = \mathcal{G}^0(\int_0^\cdot \nu(s) ds)} \dkuo{\frac{1}{2}\int_0^\cdot |\nu(s)|^2 ds}
$$
$\mathcal{G}^0$ is the solution map of \eqref{eq:4_sol}
\begin{equation}\label{eq:4_sol}
  \phi(t) = x_0 + \int_0^t \sigma(s,\phi_s)ds
\end{equation}
\end{thm}

For functional of $X(\varepsilon t )$ generally, we have the following result
\begin{thm}{}{}
Let $f \in \mathcal{C}^1_b(\mathbb{R}^d; \mathbb{R}^m)$. Then the process $f(X(\varepsilon t ))$ satisfies a LDP as $\varepsilon \rightarrow 0$ with rate function $J^f$ and speed $\varepsilon $.
\begin{equation}
  J^f(g) = \inf_{\dkuo{Df(x_0)\varphi = g}}{J(\varphi)}
\end{equation}

\end{thm}
\begin{pf}
	This proof bases on Theorem \ref{thm:LDP_PSDE} and the delta method of large deviation\citep{Gao2011DELTAMI}.
	
For any $f \in C_b^1\left(\mathbb{R}^d ; \mathbb{R}^m\right)$, define $\Phi: C_0^1\left([0, T], \mathbb{R}^d\right) \rightarrow C_0^1\left([0, T], \mathbb{R}^m\right)$ as follows:
$$
f(\varphi)(t)=f(\varphi(t)), \quad t \in[0, T] .
$$
 $\Phi$ is Hadamard differentiable and its Hadamard differential at constant function $\varphi \equiv$ $f\left(x_0\right)$ is
$$
\Phi_{f\left(x_0\right)}^{\prime}(\psi)=(D f)\left(x_0\right) \psi, \quad \psi \in C_0^\alpha\left([0, T], \mathbb{R}^d\right)
$$
Then the result follows from the delta method.
\end{pf}
\bibliographystyle{elsarticle-harv} 
\bibliography{elsarticle-template-harv}





\end{document}